\documentclass[letterpaper, 10 pt, conference]{ieeeconf}  
\pdfoutput=1

\IEEEoverridecommandlockouts                              

\usepackage{amsmath,amssymb,enumerate,amsfonts,mathrsfs,theorem}
\usepackage{graphicx}
\usepackage{color}
\usepackage[T1]{fontenc}
\usepackage{subfigure}

\theorembodyfont{\it}
\newtheorem{theorem}{Theorem}
\newtheorem{proposition}{Proposition}

\theorembodyfont{\rm}

\newtheorem{definition}{Definition}

\title{\LARGE \bf
Realization of nonlinear behaviors from organizing centers*
}

\author{Alessio Franci$^{1}$ and Rodolphe Sepulchre$^{2}$%
\thanks{*This paper presents research results of the Belgian Network DYSCO
(Dynamical Systems, Control, and Optimization),
funded by the Interuniversity Attraction Poles Programme,
initiated by the Belgian State, Science Policy Office.
The scientific responsibility rests with its authors.}
\thanks{$^{1}$Alessio Franci is with the Engineering Department,
Cambridge University, United Kingdom
        {\tt\small af529@cam.ac.uk}}%
\thanks{$^{2}$Rodolphe Sepulchre is with the Engineering Department,
Cambridge University, United Kingdom
        {\tt\small r.sepulchre@eng.cam.ac.uk}}%
}

\begin{document}

\fontencoding{T1}\fontsize{10}{10.8pt}\selectfont

\maketitle
\thispagestyle{empty}
\pagestyle{empty}

\begin{abstract}
Borrowing the concept of {\it organizing center} from singularity theory, the paper proposes a methodology to realize nonlinear behaviors such as switches, relaxation oscillators, or bursters from core circuits that reveal the fundamental role of monotonicity and feedback in their robustness and modulation.
\end{abstract}

\section{INTRODUCTION}
\label{SEC: intro}
The minimal (internal) realization of a prescribed (external) behavior is an important question 
of system theory. Beyond the realization theory of linear-time invariant systems, notable examples include 
the realization of passive systems with passive elements \cite{Willems1976}, the existence of positive state-space realizations
for  positive systems \cite{Benvenuti2004}, and the realization of Hamiltonian systems \cite{VanderSchaft1987}.

The present paper explores the realization of nonlinear behaviors organized by a singular static behavior.

As a first step of our construction,  we borrow from singularity theory  \cite{Golubitsky1985}  the concept of {\it organizing center} and the property that any robust static behavior lives in the {\it universal unfolding} of a singularity, that is, admits a normal form characterized by local algebraic conditions that identify its organizing center. In other words, any static behavior that is robust to small perturbations can be represented as  a local deformation of a singular behavior in the space of unfolding parameters. The number of unfolding (or modulation) parameters is the codimension of the singularity.

As a second step of our construction, we observe that the algebraic characterization of a given singularity and its unfolding can be translated in a feedback interconnection of monotone nonlinearities. The observation is mathematically straightforward but shows the relevance of singularity theory for behavioral theory: any static behavior that is robust to small perturbations can be realized as an interconnection of  monotone nonlinearities with interconnection laws parametrized by the unfolding parameters. The circuit representation of static normal forms reveals the key role of monotonicity and feedback in organizing static behaviors.

As a third and final step of the construction, we turn the static behaviors organized by singularities into nonlinear dynamic behaviors whose limit sets are shaped by the static behavior. The hierarchy of singularity theory between state, bifurcation parameter, and unfolding parameters translates into a hierarchy of time-scales while singular perturbation theory proves the existence of hyperbolic attractors in the vicinity of the singular behavior. 

For the sake of illustration, the present paper limits the exposition of this general procedure to the simplest situation:  (i) we only consider two of the simplest singularities (the co-dimension one hysteresis and the co-dimension three winged cusp);  (ii) we only consider the addition of first-order linear filters when moving from static to dynamic behaviors; and (iii), we rely on singular perturbation theory to infer the global analysis of the dynamical behaviors from the static analysis of the singular circuit, resulting in multiple time-scales nonlinear behaviors.

We show that this restricted framework is sufficient to recover  the most popular nonlinear behaviors: bistable switches and relaxation oscillators organized by the hysteresis singularity; rest-spike bistable models and bursters organized by the winged-cusp. The motivation for the proposed construction indeed directly comes from the recent work \cite{Franci2013b} which shows that (possibly high-dimensional) biophysical models of neurons are organized by a winged-cusp and that its circuit representation  is in one-to-one correspondence with the physiological regulatory feedback parameters that modulate bursters across a vast range of signaling behaviors.

\section{SINGULARITIES ORGANIZE STATIC BEHAVIORS}

We refer the reader to the book \cite{Golubitsky1985} for a comprehensive introduction to the singularity approach to bifurcation problems but revisit basic elements of the theory from a systems viewpoint. The theory in \cite{Golubitsky1985} is a robust bifurcation theory for scalar mappings
\begin{equation}\label{EQ: generic bif problem}
g(y,u)=0,\quad y,u\in\mathbb R,
\end{equation}
where $g$ is a smooth function. The set of pairs $(y,u)$ satisfying (\ref{EQ: generic bif problem}), that is, the graph of $g$,  defines a static behavior. The roles of $u$ and $y$ are differentiated: the variable $y$ is considered as an output variable whereas the bifurcation parameter $u$ is considered as an input (control) variable.

{\it Singular points} satisfy $g(y^\star,u^\star)=g_y(y^\star,u^\star)=0$, where $g_y:=\frac{\partial g}{\partial y}$ and similarly for higher order derivatives. Indeed, if $g_y(y^\star,u^\star)\neq0$, then the implicit function theorem applies and the graph is necessarily regular at $(y^\star,u^\star)$, therefore excluding intrinsically nonlinear static behaviors like co-existence of multiple solutions, or solution branching. Generically, at a singular point, $g_{yy},g_u\neq0$, in which case the singularity is a (codimension-$0$) fold. Locally near a fold, the graph ``looks like" the solution set of $\pm y^2+u=0$, whose left hand side defines the {\it normal form} of the fold. Higher-codimension singularities are inherently non robust to perturbations, which motivates the following definition.

\begin{definition}{\bf[Equivalence and Strong Equivalence]}\label{DEF: equivalence}\\
Two smooth mappings $h,g:\mathbb{R}\times\mathbb{R}\to\mathbb{R}$ are  {\it equivalent} near the origin if there exist smooth functions $Y(y,u)$, $U(u)$, and $Q(u,y)$, with $Y(0,0)=U(0)=0$ and $Y_y(y,u),U_u(u),Q(y,u)>0$, such that
\begin{equation}\label{EQ: equiv eq}
h(y,u)=Q(y,u)g(Y(y,u),U(u))
\end{equation}
near the origin. If (\ref{EQ: equiv eq}) holds with $U(u)=u$, the two mapping are said to be {\it strongly equivalent}.
\end{definition}

Two equivalent mappings are diffeomorphic through the diffeomorphism $(Y,U)$. The fact that $U$ does not depend on $y$ reflects the distinction between control and output: the equivalence does not allow feedback transformations at the input. Definition 1 is the key difference between the singularity approach to bifurcation problems and the classical approach based on structural stability (see \cite[Chapter 3]{Guckenheimer2002} for an excellent introduction).

The {\it recognition problem} \cite[Chapter II]{Golubitsky1985} provides a finite list of degenerate conditions (that is, higher-order zero terms in Taylor expansion beside $g=g_y=0$), that characterize a given singularity up to the equivalence relation of Definition  \ref{DEF: equivalence}. The number of degenerate conditions is called the {\it codimension} of the singularity. The higher the codimension, the more complex, or degenerate, the singularity.

The key observation of singularity theory is that a singularity restricts the persistent behaviors that can be obtained by small perturbations of the singular behavior. This is formalized through the concept of universal unfolding.

\begin{definition}{\bf[Universal unfolding]}\label{DEF: univer unfolding}
Given a smooth mapping $g(y,u)=0$ possessing a codimension $k$ singularity at the origin, a universal unfolding is a $k$-parameter family of smooth mapping $G(\cdot,\cdot;\alpha)$, $\alpha\in\mathbb R^k$, such that
\begin{itemize}
\item[]{\it a)} $G(y,u;0)=g(y,u)$;
\item[]{\it b)} for all smooth function $p(y,u)$ and sufficiently small $c>0$, there exist smooth mappings $Q,Y,U$ and $A$ such that
\end{itemize}
\begin{equation}\label{EQ: univ unf equation}
g(y,u)+c p(y,u)=Q(y,u,c)G(Y(y,u,c),U(u,c),A(c))
\end{equation}
where $Q(y,u,0)\equiv 1$, $Y(y,u,0)\equiv y$, $U(u,0)\equiv u$, and $A(0)=0$. The $k$ parameters $\alpha_1,\ldots,\alpha_k$ are called the {\it unfolding parameters}.
\end{definition}

\begin{figure}[htb]
\center
\includegraphics[width=0.45\textwidth]{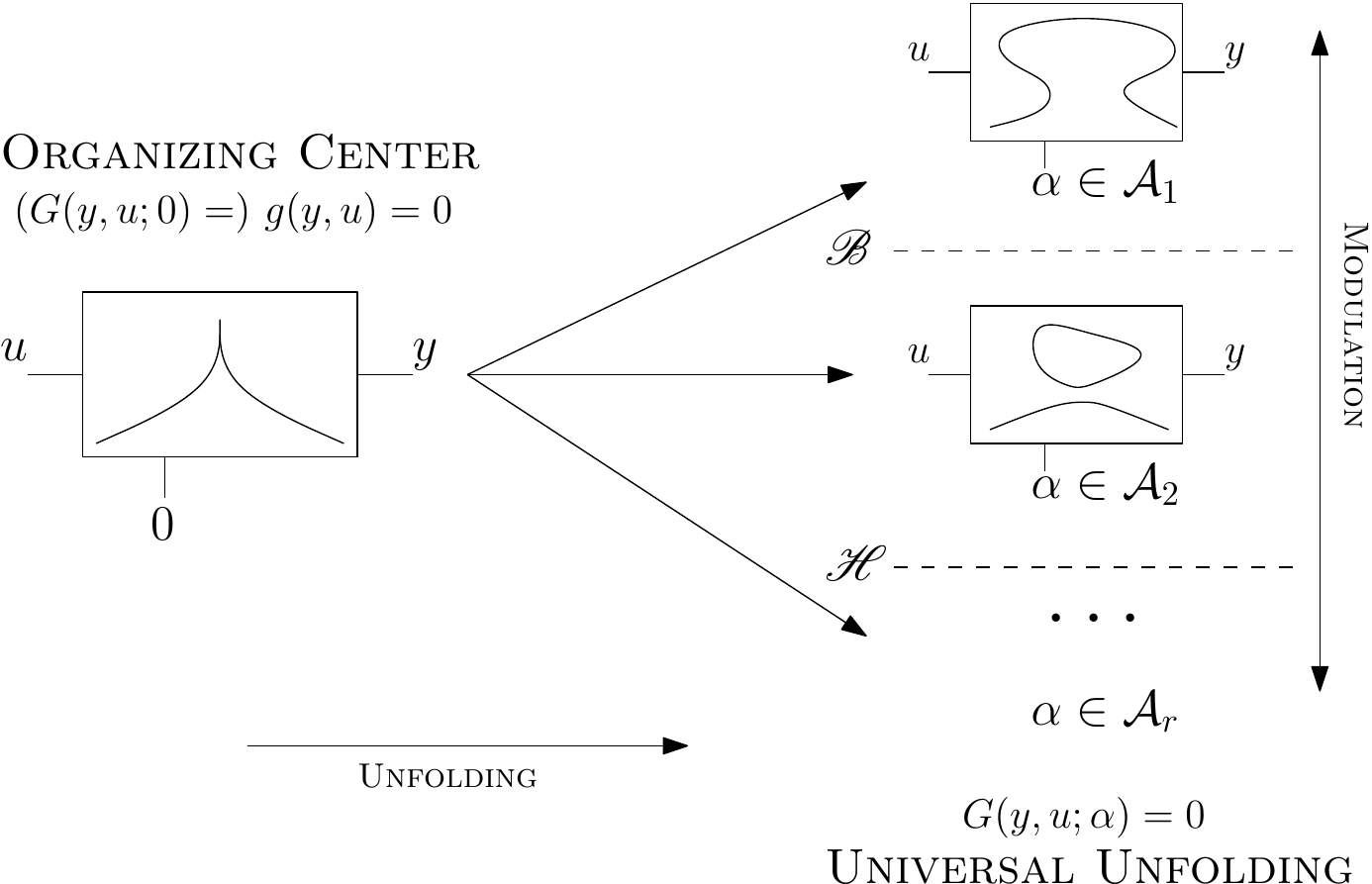}
\caption{Singularity theory characterizes all possible static behaviors in the neighborhood of a singular behavior. Distinct regions in the unfolding parameter space are separated by transition varieties, which define equivalence classes of persistent bifurcation diagrams.  The figure illustrates part of the universal unfolding of the winged cusp $-y^3-u^2$. $\mathscr B$: transcritical bifurcation transition variety. $\mathscr H$: hysteresis transition variety. $\mathcal A_i$, $i=1,\ldots,r$, denotes disconnected regions of the unfolding space, separated by transition variety, containing equivalence classes of persistent bifurcation diagrams.}\label{FIG: syst sing}
\end{figure}

The codimension of the singularity determines the number of unfolding parameters. See \cite[Section III.2]{Golubitsky1985} for an algebraic-geometric interpretation of this fact. A universal unfolding provides a {\it sensitivity chart} of the static behavior. The input $u$ is a bifurcation parameter.  In contrast, unfolding parameters are not. Rather, they shape the $u$-controlled behavior according to the universal unfolding $G$. This hierarchy is analytically reflected in the dependendence of the smooth mappings $(Y,U,A)$ on the output $y$, the input $u$, and the parameter $c$ in  (\ref{EQ: univ unf equation}). Equivalent behaviors in the sense of Definition \ref{DEF: equivalence}   are separated by {\it transition varieties}, that is, high-sensitivity codimension-$1$ algebraic varieties in the unfolding space $\mathbb R^k$, on which the static behavior possesses codimension-1 singularities.

The {\it recognition problem for universal unfolding} \cite[Chapter III]{Golubitsky1985} provides a list of simple algebraic conditions for a parametrized family of static behaviors $G(x,\lambda;\alpha)=0$ to be a universal unfolding of a given singularity.

The relevance of singularity theory to classify static behaviors is illustrated in Figure \ref{FIG: syst sing}. High-codimension singularities determine, or ``organize", a family of robust behaviors. In the terminology of  Ren\'e Thoms, each singularity is an  {\it organizing center}. By unfolding an organizing center, we obtain a finite family of qualitatively distinct, robust static behaviors,
separated by transition algebraic varieties. By modulating unfolding parameters across transition varieties we can reliably and robustly control the qualitative properties of the static behavior.

As a matter of illustration, the rest of the paper focuses on the two simplest singularities: the hysteresis and the winged cusp.  The {\it hysteresis} singularity  is the codimension-$1$  singularity with normal form
\begin{equation}\label{EQ: hysteresis singularity}
g_{hy}(y,u):=-y^3-u\,.
\end{equation}
The normal form of its universal unfolding  is
\begin{equation}\label{EQ: hysteresis singularity uu}
G_{hy}(y,u;\beta):=-y^3-u+\beta y\,.
\end{equation}

The {\it w-cusp} (winged cusp)  singularity is the codimension-$3$ singularity with normal form
\begin{equation}\label{EQ: wcusp singularity}
g_{wcusp}(y,u):=-y^3-u^2\,,
\end{equation}
The key difference between the two singularities is that the input-output relationship is {\it monotone} in the hysteresis whereas it is {\it non-monotone} in the winged cusp. The normal form of the universal unfolding of the w-cusp  is
\begin{equation}\label{EQ: wcusp singularity uu}
G_{wcusp}(y,u;\alpha,\beta,\gamma):=-y^3-u^2+\alpha+\beta y+\gamma u y\,.
\end{equation}

\section{FROM UNIVERSAL UNFOLDINGS TO STATIC CIRCUITS}

A straightforward but key message of the present paper is that singular behaviors and their universal unfoldings can be realized as interconnections of monotone sigmoidal nonlinearities. Converting the algebraic characterization to a circuit characterization is insightful and instrumental in identifying organizing centers  of static behaviors from interconnection structures.  For a smooth function $f(u)$ and $n\geq0$, we use the notation $f^{(n)}:=\frac{d^{n} S}{d u^{n}}$.

\begin{definition}{\bf [Sigmoidal nonlinearity]}\label{DEF: smooth sat}
The {\it sigmoidal static behaviors} is defined by
$$-y+S(u)=0,$$
where $S$ is a smooth function satisfying the following properties:
\begin{itemize}
\item[{\it a.}] $S(-u)=-S(u)$, for all $u\in\mathbb R$ (odd);
\item[{\it b.}] $S^{(1)}(u)>0$, for all $u\in\mathbb R$ (monotone);
\item[{\it c.}] $\lim_{u\to\pm\infty}S^{(1)}(u)=0$ (saturated);
\item[{\it d.}] For all $n\in\mathbb N$, $S^{(2n+1)}(0)\neq 0$ with $S^{(1)}(0)=1$ (regular);
\item[{\it e.}] For all $n\in\mathbb N$, ${\rm sgn}\left(S^{(2n)}(u)\right)=-{\rm sgn}(u)$, for all $u\neq 0$.
\end{itemize}
\end{definition}
Sigmoidal nonlinearities are widespread in nonlinear models encountered in engineering, physics, and biology: they include sigmoids, saturations, and hyperbolic tangents.

\subsection{A static circuit realization of the hysteresis}
The basic observation  
\begin{equation}\label{EQ: basic 3 by sat}
y-S(y)=-\frac{S^{(3)}(0)}{6!}y^3+\mathcal O(y^4),
\end{equation}
suggests that the cubic hysteresis normal form (\ref{EQ: hysteresis singularity}) can be realized (in the strong equivalence sense) as a positive feedback loop around a sigmoidal nonlinearity, {\it i.e.}
$$-y+S(y+u)=0.$$

\begin{figure}[h!]
\center
\includegraphics[width=0.35\textwidth]{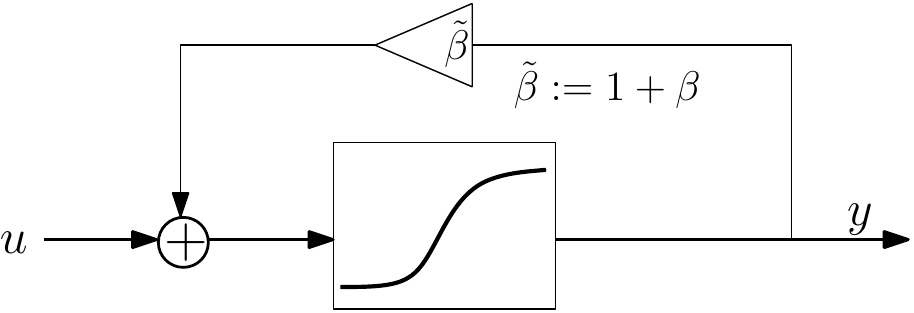}
\caption{Universal unfolding of the hysteresis by interconnection.}\label{FIG: hy int}
\end{figure}
\begin{proposition}\label{PROP: hy int}
For $\beta=0$, the static behavior $-y+S(y+u)=0$ of the feedback loop in Figure \ref{FIG: hy int} is strongly equivalent to the hysteresis singularity (\ref{EQ: hysteresis singularity}).
\end{proposition}
{\bf Proof.} We prove this proposition by solving the recognition problem for the hysteresis singularity \cite[Proposition II.9.1]{Golubitsky1985}. Let $h(y,u):=-y+S(y+u)$. We have to verify that
$$h(0,0)=h_y(0,0)=h_{yy}(0,0)=0,$$ $$ h_{yyy}(0,0)<0,\ h_{u}(0,0)<0,$$
which are all easily verified by invoking properties of the sigmoidal behavior $S$.\hfill$\square$\\

\begin{proposition}\label{PROP: hy int uu}
The static behavior $-y+S(y+u+\beta y)=0$ of the feedback loop in Figure \ref{FIG: hy int} is a universal unfolding of the hysteresis singularity, with $\beta$ as unfolding parameter. Its persistent static behaviors are given by the persistent bifurcation diagrams of the hysteresis \cite[Page 205]{Golubitsky1985}.
\end{proposition}
{\bf Proof.} We apply the recognition problem for the universal unfolding of the hysteresis \cite[Proposition 4.4]{Golubitsky1985}. Let $H(y,u;\beta):=-y+S(y+u+\beta y)$. $H$ is a one parameter unfolding of $h$, defined in Proposition \ref{PROP: hy int}, and, by  Proposition \ref{PROP: hy int}, $h$ is strongly equivalent to the hysteresis. It remains to prove that
$
\text{det}\left(\begin{array}{cc}
h_{u} & h_{u y}\\ H_{\beta} & H_{\beta y}
\end{array}\right)\neq0\,,
$
at $y=u=\beta=0$, which easily follows by the sigmoidal properties of $S$.
\hfill$\square$\\

The elementary proofs illustrate the relevance of a theory robustly characterizes nonlinear behaviors from local calculations.

\subsection{The winged-cusp static behavior by interconnection}
The non-monotone behavior $u^2$ characterizing the w-cusp singularity (\ref{EQ: wcusp singularity}) and its universal unfolding (\ref{EQ: wcusp singularity uu}) is obtained by cascading the realization of the hysteresis in Figure \ref{FIG: hy int} with a non-monotone static nonlinearity. 

\begin{figure}[h!]
\center
\includegraphics[width=0.40\textwidth]{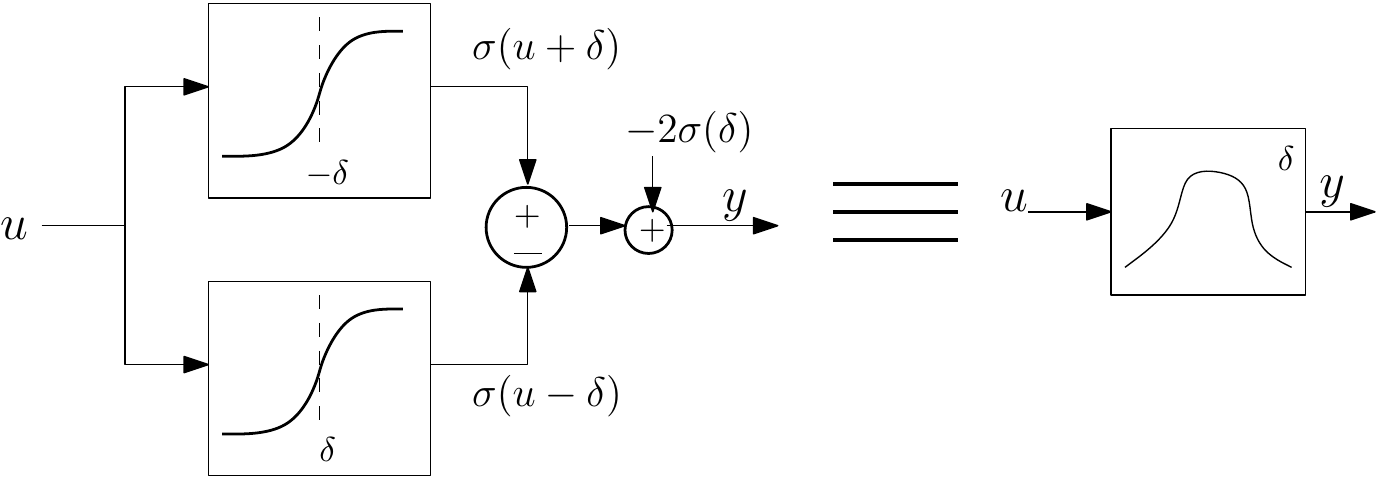}
\caption{Parallel interconnection of two monotone behavior realizes a non-monotone behavior.}\label{FIG: basic par 2}
\end{figure}

In essence, a non-monotone static behavior  can be realized as the negative parallel interconnection of two monotone behaviors. We define the {\it bump} nonlinearity
\begin{equation}\label{EQ: bump def}
B_\delta(u):=S(u+\delta)-S(u-\delta)-2S(\delta),\quad \delta\neq0,
\end{equation}
sketched in Figure \ref{FIG: basic par 2}.

\begin{figure}[h!]
\center
\includegraphics[width=0.40\textwidth]{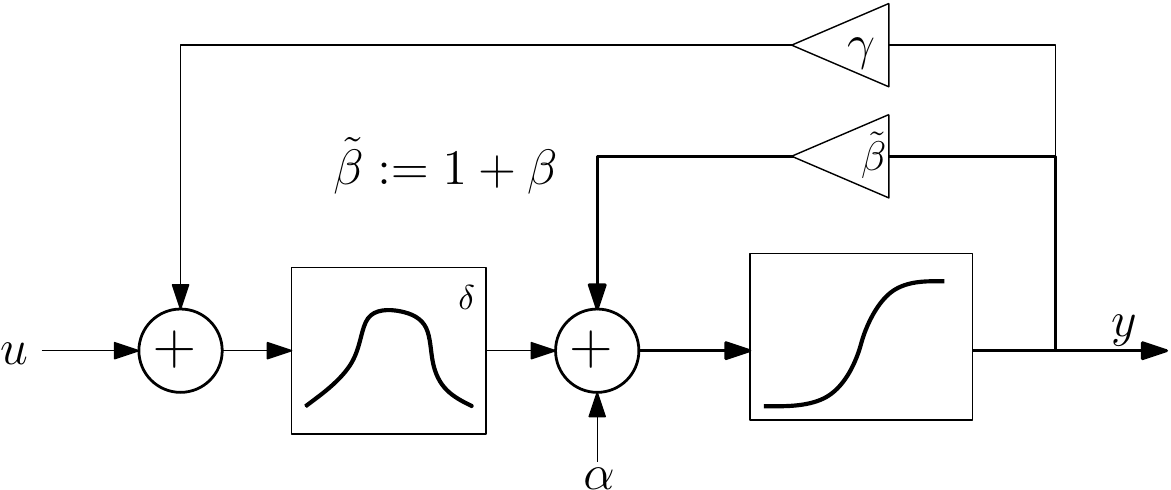}
\caption{Winged cusp singularity by interconnection.}\label{FIG: wcusp int}
\end{figure}

From properties {\it d.} and {\it e.} of $S$ in Definition \ref{DEF: smooth sat}, it holds that, for $\delta$ sufficiently close to the origin,
\begin{align*}
B_\delta(0)=0,\ B_\delta^{(1)}(0)=0,\ B_\delta^{(2)}(0)=2S^{(2)}(\delta)\neq0,
\end{align*}
and therefore
$$B_\delta(u)=2S^{(2)}(\delta)u^2+\mathcal O(u^4).$$
By cascading $B_\delta$ with the circuit in Figure \ref{FIG: hy int}, the resulting static behavior reads
$$-y+S(B_\delta(u)+y+\beta y)=0, $$
which is a good candidate to realize the universal unfolding of the w-cusp. The additive unfolding parameter $\alpha$ in (\ref{EQ: wcusp singularity uu}) can then just be added to the input of the sigmoid. The unfolding parameter $\gamma$ involves multiplication of the output and input. The equation
$$B_\delta(u+\gamma y/2)=2S^{(2)}(\delta)(u^2+\gamma u y+\gamma^2 y^2/4)+\mathcal O(3),$$ 
shows that it can be realized by a feedback loop around the bump nonlinearity.  By the recognition problem for the universal unfolding of the winged cusp, the term $\gamma^2 y^2/4$ and higher order terms $\mathcal O(3)$ will play no qualitative role in the unfolding. The overall interconnection structure is summarized in Figure \ref{FIG: wcusp int} and we have the following results, whose proofs easily follow from the recognition problem for the winged-cusp and its universal unfolding.

\begin{proposition}\label{PROP: wcusp int}
For $\delta>0$ and sufficiently small, the static behavior $-y+S(B_\delta(u)+y)=0$ of the interconnection in Figure \ref{FIG: wcusp int} with $\alpha=\beta=\gamma=0$ is strongly equivalent to the winged-cusp (\ref{EQ: wcusp singularity}).
\end{proposition}
\begin{proposition}\label{PROP: wcusp int uu}
For $\delta>0$ and sufficiently small, the static behavior $-y+S(B_\delta(u+\gamma y/2)+y +\alpha+\beta y)=0$ of the interconnection in Figure \ref{FIG: wcusp int} is a universal unfolding of the winged cusp singularity with $\alpha,\beta,\gamma$ as unfolding parameters. Its persistent static behaviors are given by the persistent bifurcation diagrams of the winged cusp \cite[Page 208]{Golubitsky1985}.
\end{proposition}

\section{FROM STATIC TO DYNAMIC BEHAVIORS}
We explore the dynamic behaviors that are obtained by the addition of simple first-order filters in the static feedback interconnections. Most of the qualitative results will however persist when those simple linear systems are replaced by arbitrary monotone systems \cite{Angeli2003}.

The hierarchy of singularity theory between state (or output) $y$, bifurcation parameter (or input) $u$, and unfolding (or modulation) parameters $\alpha, \beta, \dots$ dictates a corresponding hierarchy of time scales. When a bifurcation parameter is varied slowly in a dynamical system, the bifurcation analysis becomes the quasi-steady state analysis of a two time-scale behavior in which the dynamics of the state is fast and the dynamics of the bifurcation parameter is slow. Ultra-slow modulation of an unfolding parameter introduces a third time-scale in which the attractors of the slow-fast system become modulated across the unfolding parameter space. Singular perturbation theory provides the necessary framework to  ensure that the singular structure captured by the static interconnection persists when the different time scales are sufficiently separated. By monotonicity of first order filters, the obtained multi-scale dynamic circuit and its critical manifold are organized by the same singularity as the underlying static behavior, which allows to analyze and tune the dynamical behavior on the same unfolding space.

The remainder of this section is provided without proofs but  mimics the construction proposed in our recent paper \cite{Franci2013b}, where the reader is referred for detailed proofs.

\subsection{Dynamical behaviors organized by the hysteresis: bistability, relaxation oscillations, and excitability}

\subsubsection{Bistability}
Adding a first-order filter with transfer function $H_f(s):=\frac{1}{\varepsilon_f s+1}$, $0<\varepsilon_f\ll 1$, in the algebraic loop in Figure \ref{FIG: hy int}, results in the simplest dynamical model of bistability with state-space realization
\begin{IEEEeqnarray}{rCl}\label{EQ: hyst fast}
\dot x &=&-x+S(x-u+\beta x)\IEEEyessubnumber\\
y&=&x \IEEEyessubnumber
\end{IEEEeqnarray}
which is reminiscent of well-studied behaviors of autocatalysis or switching behaviors, e.g. the Hopfield neuronal model \cite{Hopfield1982}, the fast behavior of Hodgkin-Huxley model \cite{Hodgkin1952},  or the genetic switch model of Griffith \cite{Griffith1968,Griffith1968a}.  

\begin{figure}[h!]
\subfigure{
\includegraphics[width=0.35\textwidth]{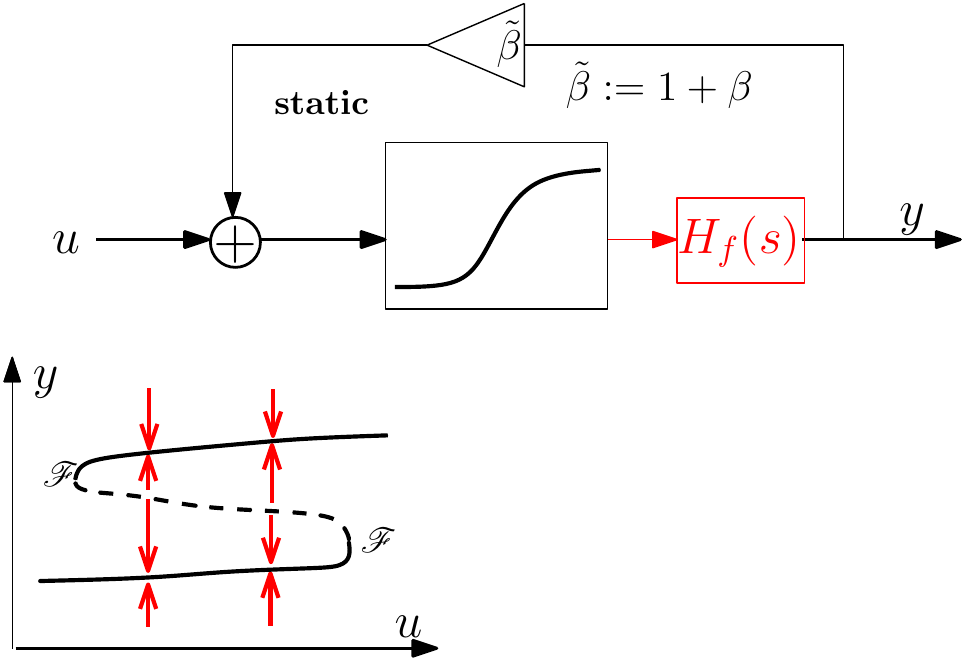}
}\\
\vspace{-3cm}
\flushright{
\subfigure{
\includegraphics[width=0.23\textwidth]{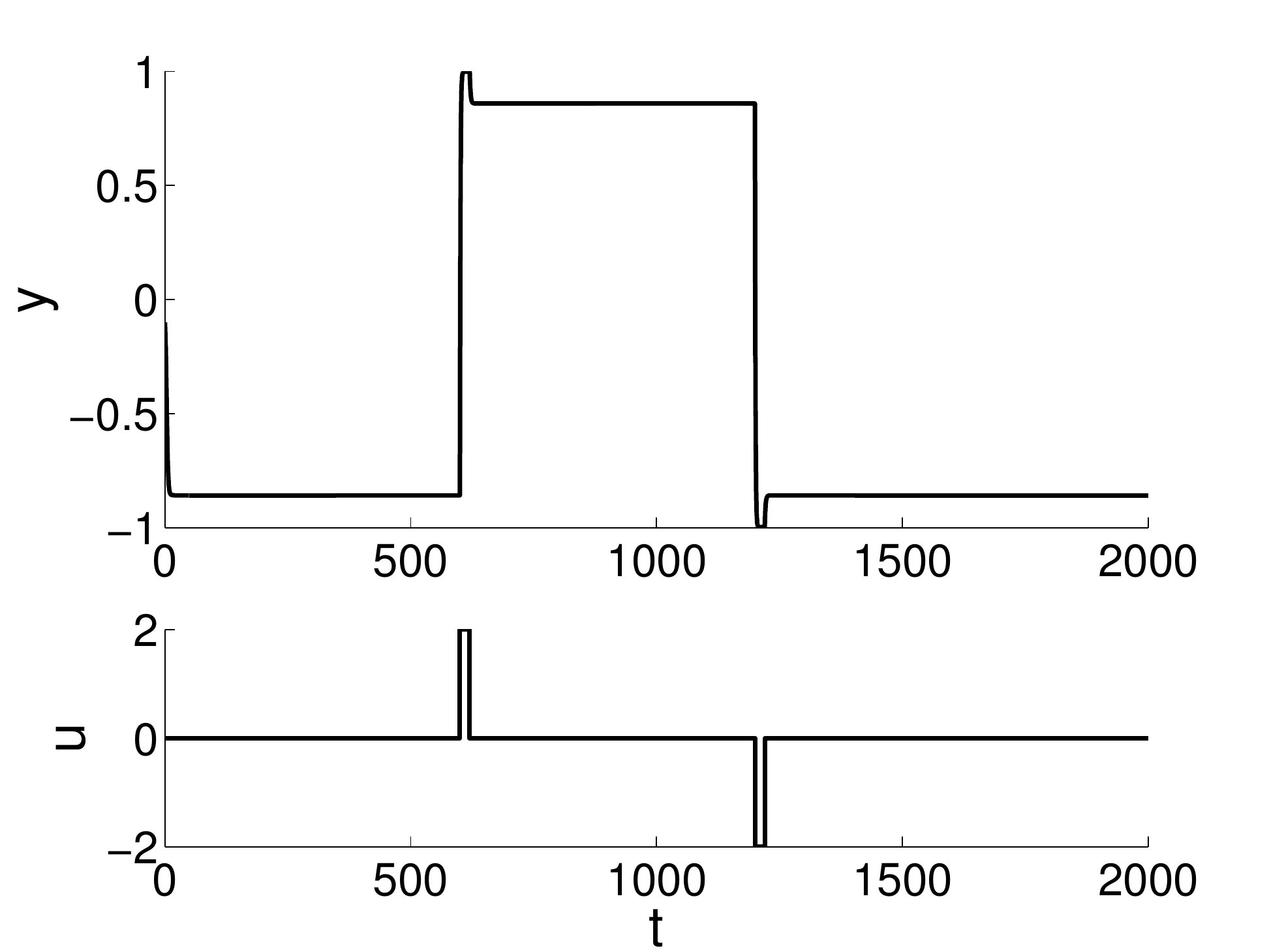}
}
}
\caption{{\bf Bistability in a circuit organized by the hysteresis}. Top: circuit realization. Bottom left: sketch of the dynamical circuit behavior for $\beta>0$. Branches of stable steady state are drawn as full lines, whereas branches of unstable steady states as dashed lines. $\mathscr{F}$ denotes a fold singularity. Bottom right: pulse response in the bistable regime, for $\beta=0.5$. Here and in all subsequent simulations, $S(\cdot)=\tanh(\cdot)$.}\label{FIG: hy fast}
\end{figure}

The possible dynamical behaviors of (\ref{EQ: hyst fast}) are predicted from the static behavior:
\begin{itemize}
\item $\beta<0$: monostability (persistent);
\item $\beta=0$: {hysteresis bifurcation};
\item $\beta>0$: bistability (persistent - see Figure \ref{FIG: hy fast} bottom);
\end{itemize}
The same conclusions apply if the linear filter is replaced by an arbitrary monotone system, see e.g. \cite{Angeli2004}.

\subsubsection{Relaxation oscillations and excitability}
\label{SSSEC: hy slow-fast}
\begin{figure}[h!]
\subfigure{
\includegraphics[width=0.35\textwidth]{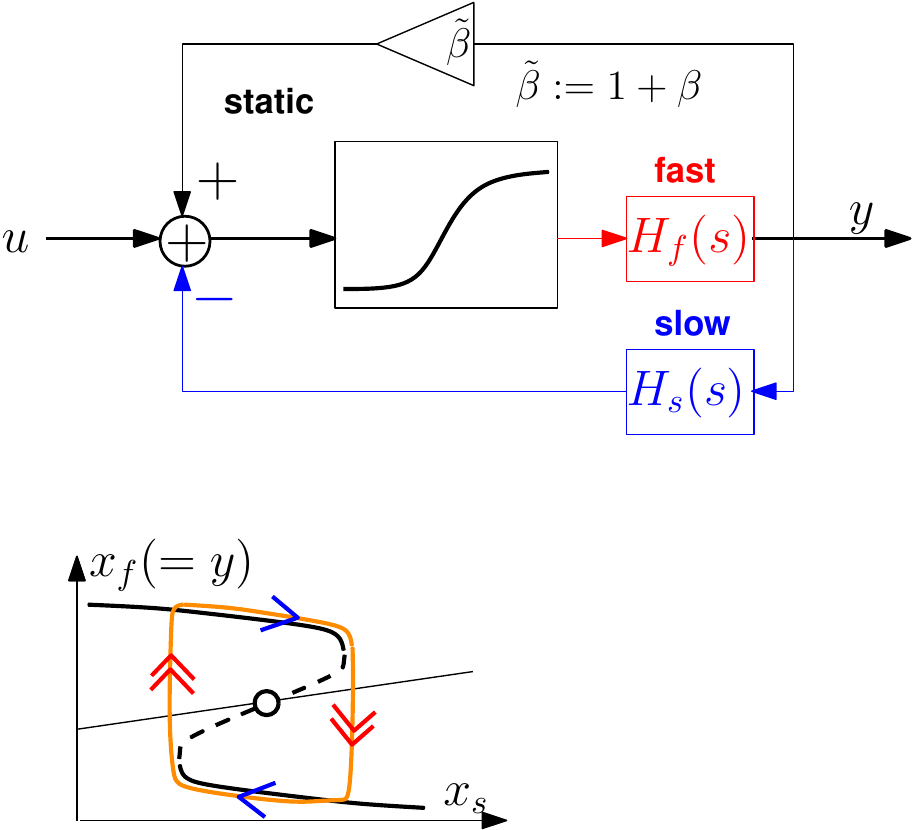}
}

\vspace{-2.75cm}
\flushright{
\subfigure{
\includegraphics[width=0.24\textwidth]{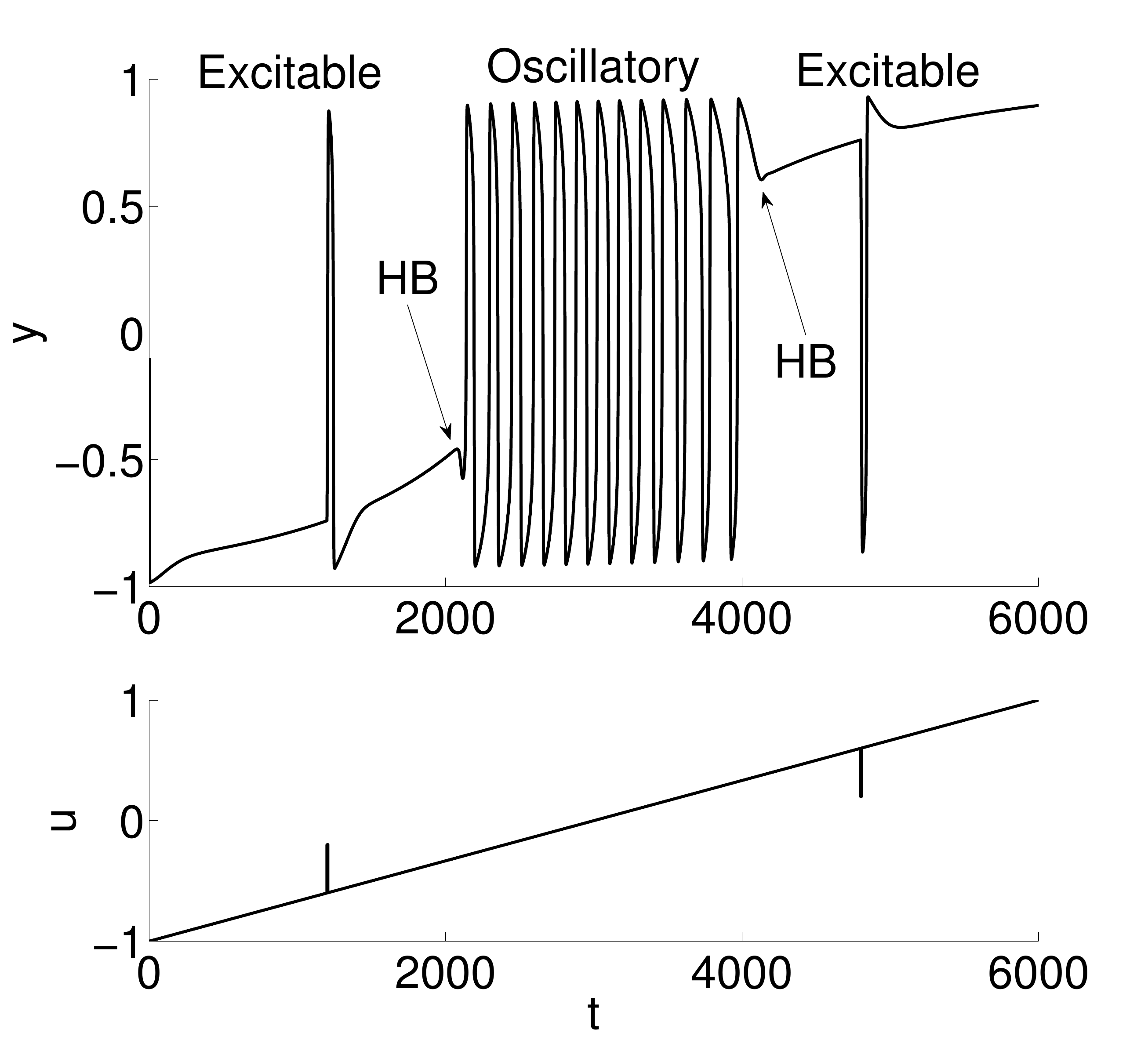}
}}
\caption{{\bf Relaxation oscillations and bistability in a slow-fast system organized by the hysteresis}. Top: circuit realization. Bottom left: sketch of the underlying slow fast phase-plane for $0<\beta<1$ and $u\sim 0$. Stable branches of the critical manifold $\dot x_f=0$ are drawn as solid thick lines, unstable ones as dashed thick lines. The slow nullcline $\dot x_s=0$ is drawn as thin line. The unstable fixed point is drawn as a circle. The orange line sketches the oscillation limit cycle, composed of slow (blue arrows) and fast (double red arrows) portions. Bottom right: circuit response for $\beta=0.5$ and $\varepsilon_f=0.01$ to an input consisting of a ramp plus two small pulses during the excitable regime. Note how the small pulses induces large responses (excitability). HB denotes a Hopf bifurcation.}\label{FIG: hy slow-fast}
\end{figure}

A bistable switch is known to be converted into a relaxation oscillator by "slow adaptation". To this end, we add a slow negative feedback loop with unitary gain modulating the control parameter (see Figure \ref{FIG: hy slow-fast} top)
\begin{IEEEeqnarray}{rCl}\label{EQ: hyst slow fast}
\varepsilon_f\dot x_f&=&-x_f+S(x_f-(u+x_s)+\beta x_f)\IEEEyessubnumber\\
\dot x_s&=&x_f-x_s\IEEEyessubnumber\\
y&=&x_f\IEEEyessubnumber
\end{IEEEeqnarray}
obtaining a fast-autocatalytic loop modulated by a slow negative feedback loop. The behavior of (\ref{EQ: hyst slow fast}) is again organized by a hysteresis singularity but, due to the slow negative feedback loop, the hysteresis singularity is now at $\beta=1$.

A global analysis of the dynamical behavior of (\ref{EQ: hyst slow fast}) as the control and the unfolding parameters are varied is straightforward near the singular limit by noticing that, by construction, its critical manifold
\begin{align*}
\mathcal M_{u,\beta}:=&\\
\{(x_f,x_s)\in\mathbb R^2:& -x_f+S(x_f+(u-x_s)+\beta x_f)=0\}
\end{align*}
is the static behavior of Proposition \ref{PROP: hy int uu} and therefore is also a universal unfolding of the hysteresis controlled by $\tilde u:=u-x_s$. In the singular limit $\varepsilon_f\to0$, the global behavior of (\ref{EQ: hyst slow fast}) is fully captured by unfolding the hysteresis singularity organizing its critical manifold. By geometric singular perturbation arguments \cite{Fenichel1971,Krupa2001a,Krupa2001b,Krupa2001c}, the same qualitative behavior persists for $\varepsilon_f>0$. For $\beta>0$ the critical manifold is multivalued and the slow negative feedback transform by hysteretic modulation the fast bistable regime in Figure \ref{FIG: hy fast} into a relaxation oscillation.

\begin{theorem}\label{THM: hy rel osc}
Consider the slow-fast dynamics (\ref{EQ: hyst slow fast}). For all $0<\beta<1$, there exists $\bar\varepsilon_f,\bar u>0$, such that, for all $|u|<\bar u$, there exists a unique exponentially unstable steady-state surrounded by an almost globally exponentially stable relaxation limit cycle.
\end{theorem}
Figure \ref{FIG: hy fast} (bottom left) sketches a phase plane construction of the relaxation limit cycle around the bistable static behavior. A rigorous proof follows along the same line as \cite{Grasman1987},\cite{Mishchenko1980},\cite{Krupa2001c}. The same result holds when first order filters are substituted by monotone systems \cite{Angeli2008}.

The dynamical behaviors of (\ref{EQ: hyst slow fast}) as $u$ varies can also easily be derived in the singular limit. In the regime $0<\beta<1$, by increasing or decreasing $u$ outside the region predicted by Theorem \ref{THM: hy rel osc} the system undergoes a Hopf bifurcation \cite[Section 3.4]{Guckenheimer2002}\cite[Section 4.7]{Chow1994} at which the unique steady state becomes exponentially stable (Figure \ref{FIG: hy fast} bottom right). The behavior is however "excitable"
\cite{FitzHugh1961,Rinzel1989,Izhikevich2007}:
the large transient response (spike) is a manifestation of the bistability of the fast subsystem. The switch is however only transient because of the adaptation provided by the slow subsystem.

\subsection{Dynamical behaviors organized by the winged-cusp: rest-spike bistability and bursting}

\subsubsection{Rest-spike bistability}
\begin{figure}[h!]
\subfigure{
\includegraphics[width=0.35\textwidth]{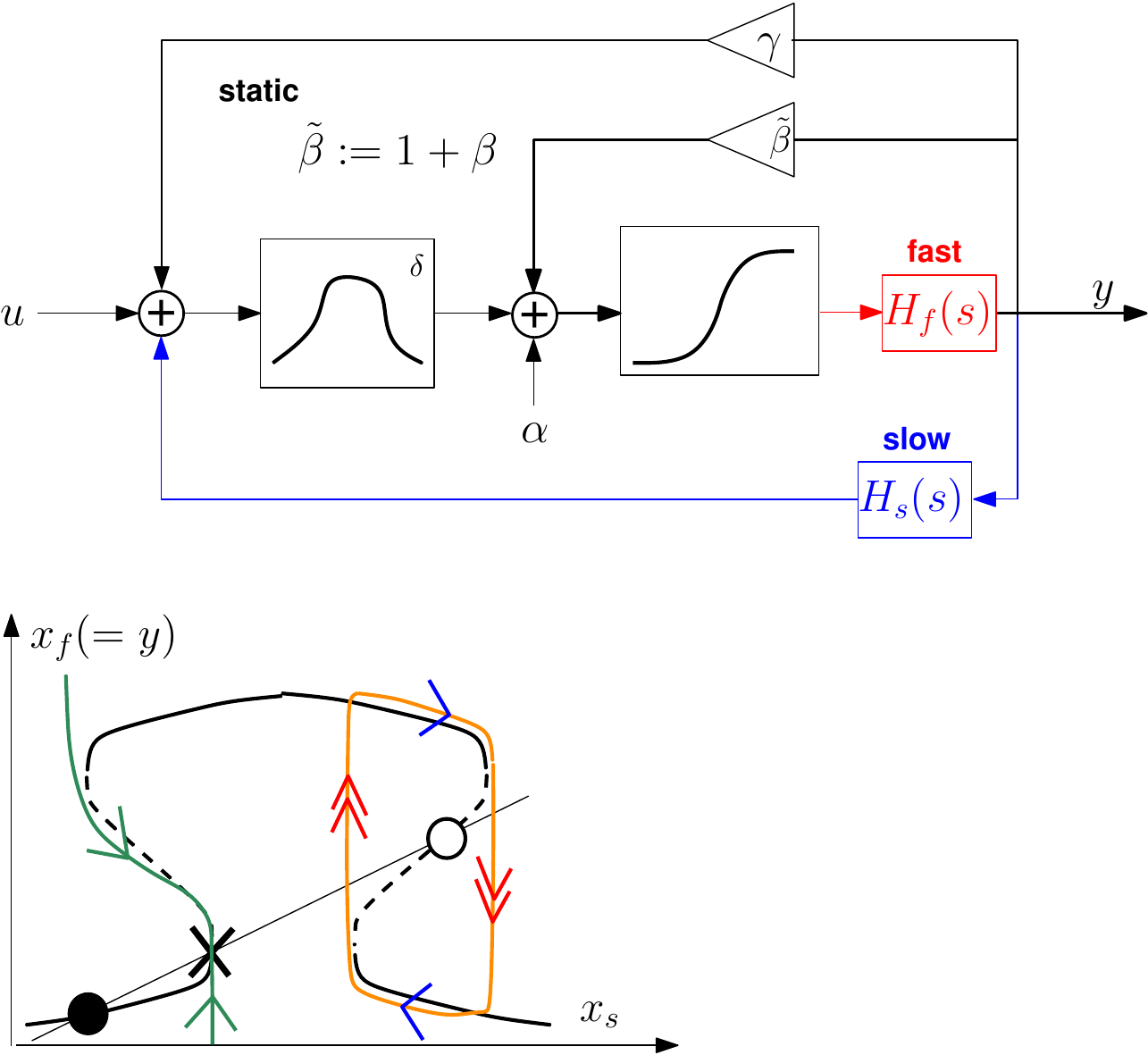}
}

\vspace{-3cm}
\flushright
\subfigure{
\includegraphics[width=0.24\textwidth]{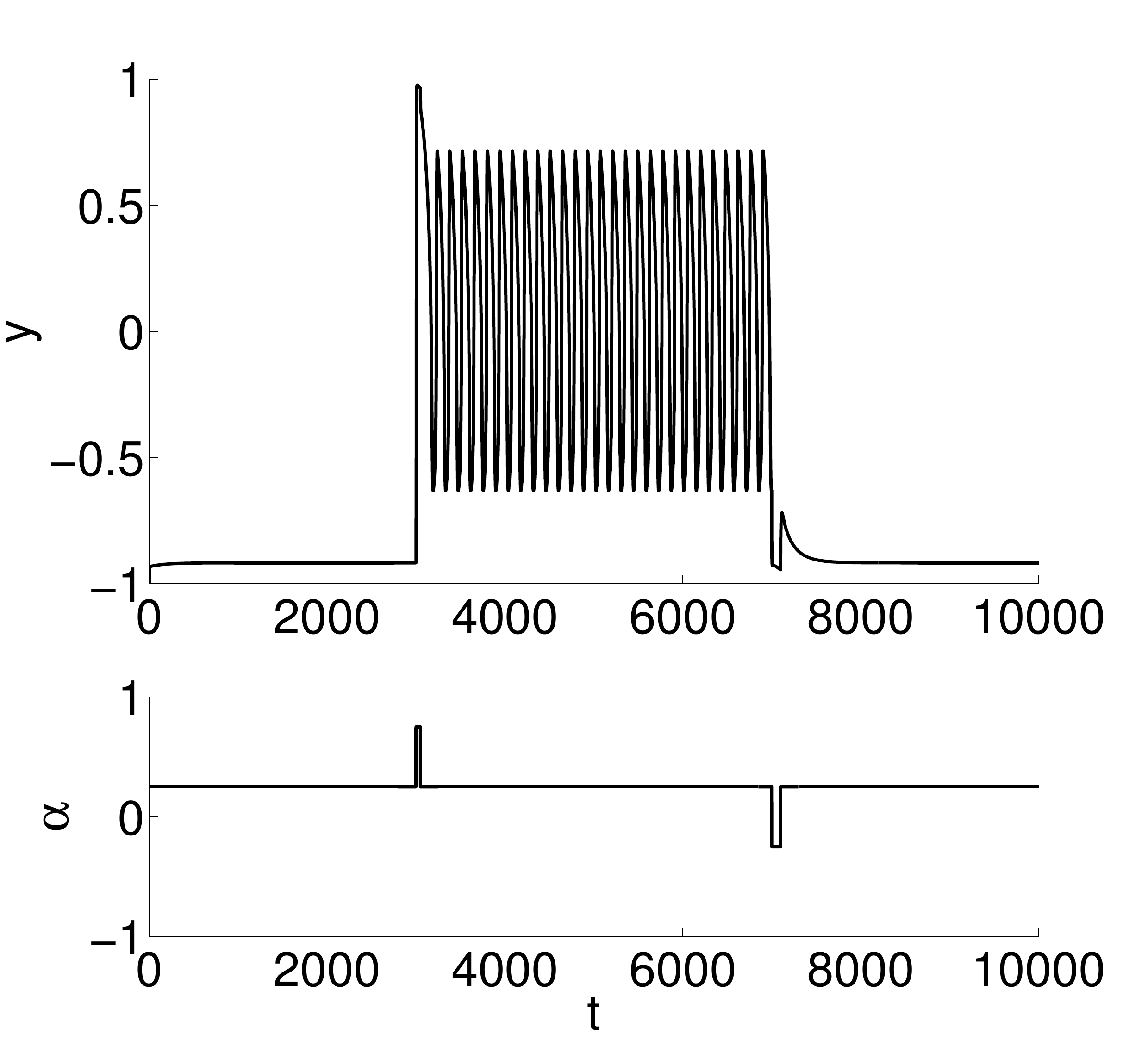}
}
\caption{{\bf Rest-spike bistability in a circuit organized by the winged-cusp.} Top: circuit realization. Bottom left: slow-fast phase-portrait organized by the mirrored hysteresis critical manifold (thick line). The this line is the slow variable nullcline. The stable manifold of the saddle point (cross) is drawn as a green oriented line. The stable steady state is drawn as full circle. Bottom right: system response to transient perturbation to the unfolding parameter $\alpha$ and $u=0.5,\,\beta=0.5,\,\gamma=1,\,\varepsilon_f=0.0075$.}\label{FIG: wcusp slow-fast}
\end{figure}

By mimicking the same construction as for the hysteresis, the slow-fast circuit organized by this singularity is depicted in Figure \ref{FIG: wcusp slow-fast} top and its state-space representation is given by
\begin{IEEEeqnarray}{rCl}\label{EQ: wcusp slow fast}
\varepsilon_f\dot x_f&=&-x_f+\IEEEyessubnumber\\
&&S\Big(x_f+B_\delta\left(u+ x_s+\frac{\gamma}{2}x_f\right)+\beta x_f+\alpha\Big)\nonumber\\
\dot x_s&=&x_f-x_s\IEEEyessubnumber\\
y&=&x_f\IEEEyessubnumber
\end{IEEEeqnarray}
We focus our analysis on an interesting dynamical behavior organized by the winged cusp, the ``rest-spike bistability'', by fixing the unfolding parameters $\beta,\gamma$ and varying $\alpha$ around the transcritical transition variety in its universal unfolding.

For suitable unfolding parameters (depending on the specific functional form of $S(\cdot)$) the static behavior is a {\it mirrored hysteresis} persistent bifurcation diagram, introduced in \cite{Franci2013b} (see Figure \ref{FIG: wcusp slow-fast} bottom left). The mirrored hysteresis organizes the rest-spike bistable behavior.
\begin{theorem}\label{THM: cusp bist}
For all $\beta>0$, there exist open sets of control ($u$) and unfolding ($\alpha,\gamma$) parameters near the pitchfork singularity in the universal unfolding of the winged cusp organizing the critical manifold of (\ref{EQ: wcusp slow fast}), in which, for sufficiently small $\varepsilon_f>0$, model (\ref{EQ: wcusp slow fast}) exhibits the coexistence of an exponentially stable fixed point and an exponentially stable relaxation limit cycle. Their basins of attraction are separated by the stable manifold of a hyperbolic saddle.
\end{theorem}
Under condition of Theorem \ref{THM: cusp bist}, short-lasting perturbations induce transitions from rest to periodic oscillations, as shown in Figure \ref{FIG: wcusp slow-fast} bottom right by transiently perturbing the unfolding parameter $\alpha$ across the transcritical bifurcation variety.

\subsubsection{Bursting}
Robust co-existence of rest and spiking is at the basis of bursting {\cite{Rinzel1987b}: an (ultra-slow) negative feedback, transforms the two-timescale bistable behavior into a three-timescale monostable behavior, much in the same way as negative feedback around the hysteresis transform a fast bistable behavior into a slow-fast oscillation.

\begin{figure}[h!]
\center
\subfigure{
\includegraphics[width=0.35\textwidth]{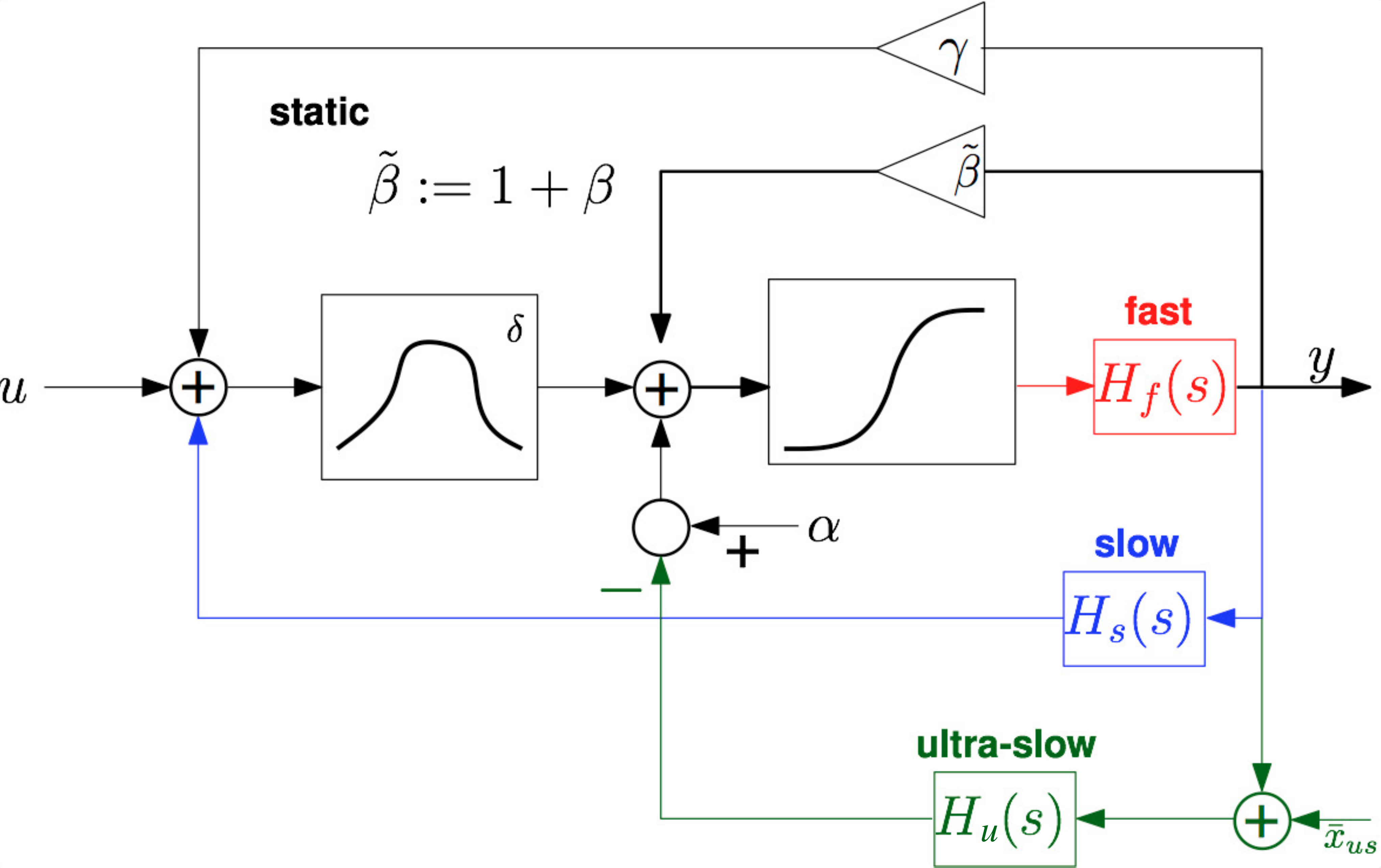}
}\\
\vspace{-1cm}
\subfigure{
\includegraphics[width=0.15\textwidth]{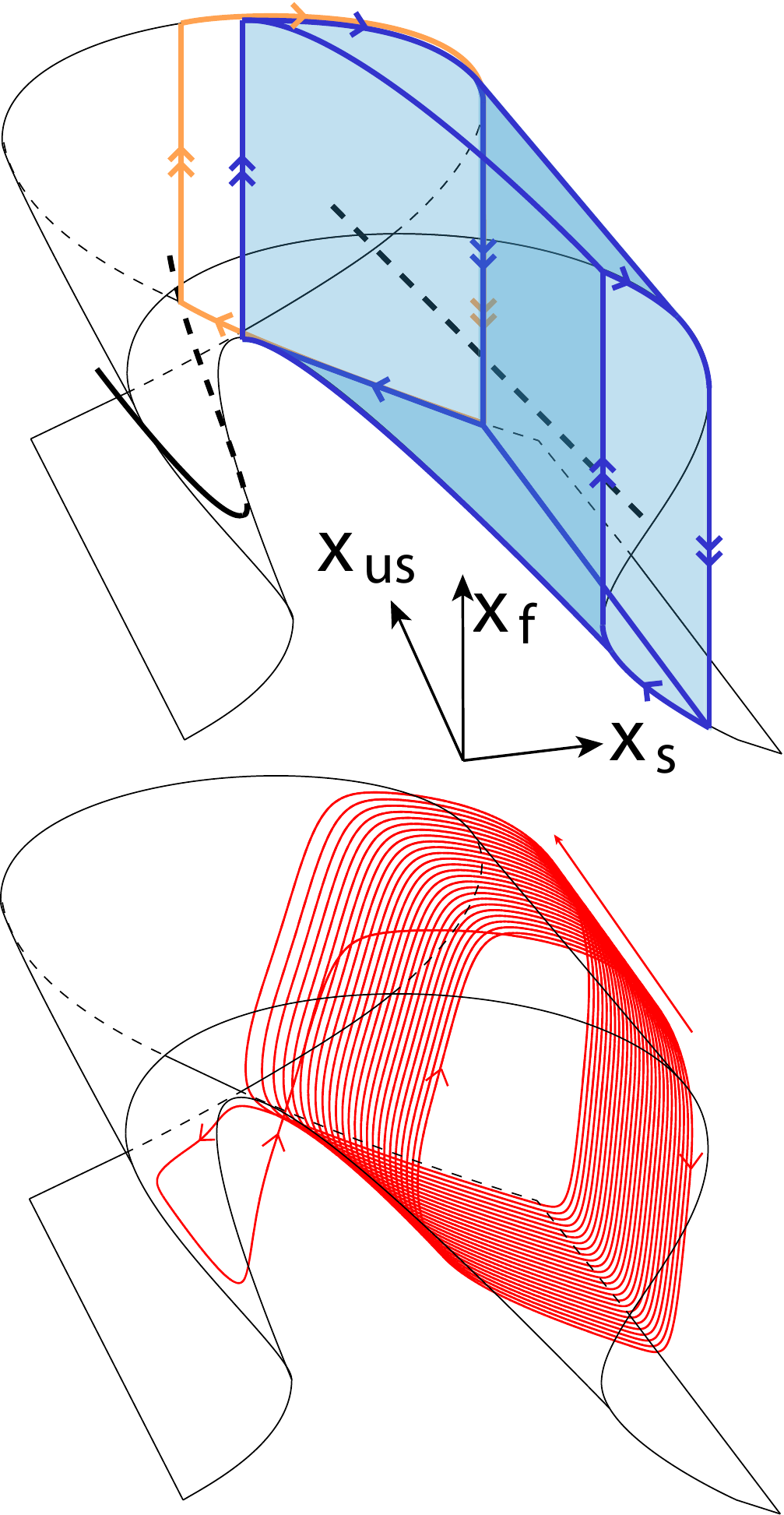}
\includegraphics[width=0.28\textwidth]{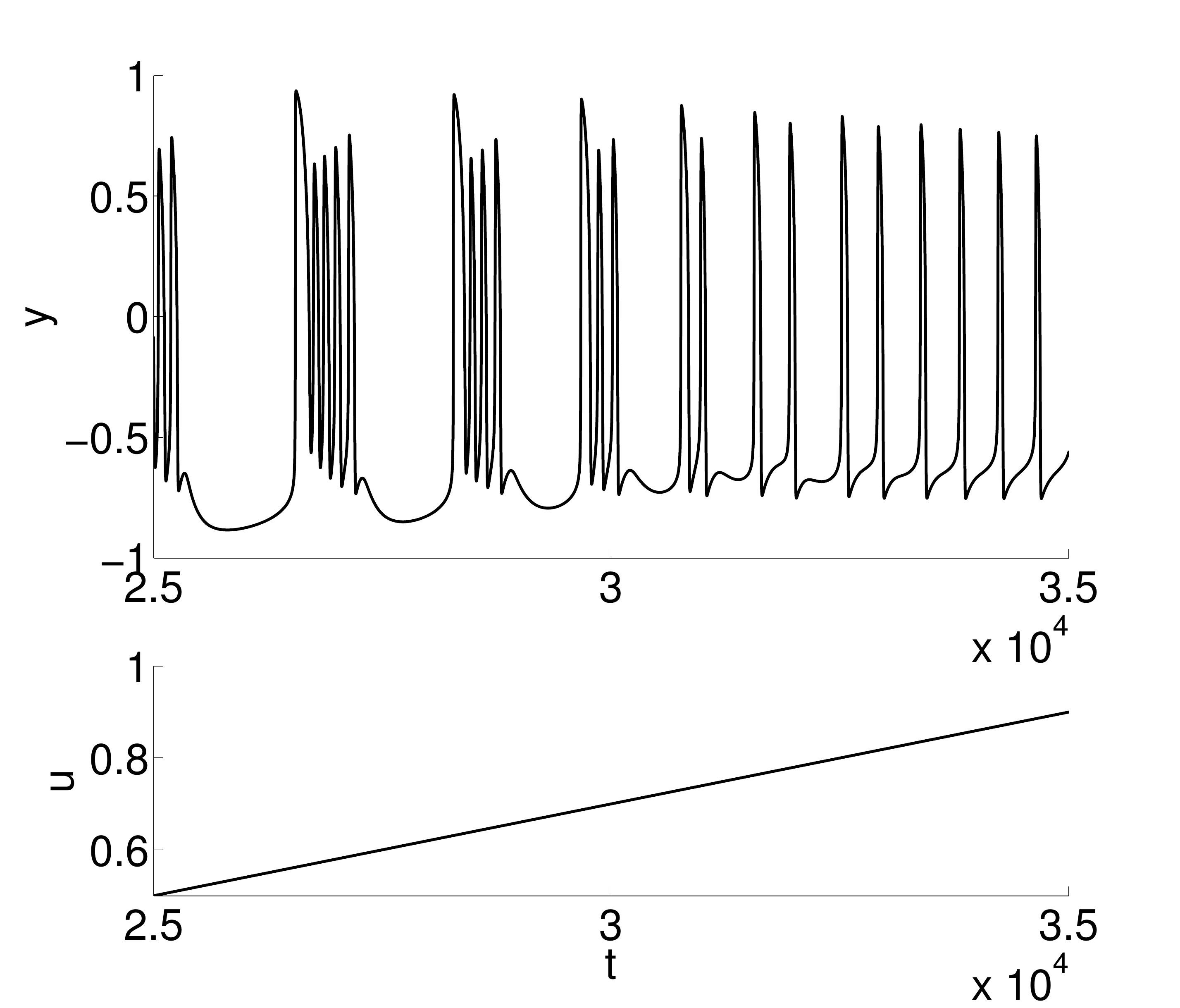}
}
\caption{A three-timescale behavior in a circuit organized by the winged cusp: bursting. Top: circuit realization. $H_u(s):=\frac{k_u}{s/\varepsilon_u+1}$. Bottom left: sketch of the three-dimensional critical manifold organizing bursting (see \cite[Figure 4]{Franci2013b} for details). Bottom right. System response to an input ramp: transition from bursting to tonic spike across a transcritical-homoclinic bifurcation. $k_u=5,\,\varepsilon_u=75^{-1},\,\bar x_u=2.5$, other parameters as in Figure \ref{FIG: wcusp slow-fast}.}\label{FIG: burst int}
\end{figure}

The unfolding parameter $\alpha$ modulates monotonically the winged-cusp singularity. It is therefore a good candidate to serve as an adaptation parameter in an ultra-slow negative feedback loop (Figure \ref{FIG: burst int} top). Effects of this ultra-slow modulation are fully captured in the unfolding space of the winged-cusp organizing center (Fig. \ref{FIG: burst int}, bottom left, and \cite[Figure 4A]{Franci2013b}).
A rigorous proof of the following theorem on a topologically equivalent normal form can be found in \cite[Theorem 2.3]{Franci2013b}.
\begin{theorem}\label{THM: cusp burst}
For all $\beta>0$, there exists an open set of control ($u$) and unfolding ($\alpha,\gamma$) parameters near the pitchfork singularity in the universal unfolding of the winged cusp organizing the critical manifold of (\ref{EQ: wcusp slow fast}) {such that, for all $\lambda,\alpha,\gamma$ in those sets, there exist $k_u,\bar x_u,\in\mathbb R$ such that,} for sufficiently small $\varepsilon_f,\varepsilon_u>0$, the circuit in Figure \ref{FIG: burst int} possesses a hyperbolic bursting attractor.
\end{theorem}
By increasing the control parameter outside the region where relaxation oscillations and rest coexist, the system behavior gradually changes from a bursting behavior to a slow spiking behavior. At the transition, the saddle-node and the homoclinic bifurcation delimiting the rest-spike bistable range merges at at transcritical-homoclinic bifurcation. See \cite[Section 3]{Franci2013b} for a geometric description of this phenomenon.

%

\section{CONCLUSIONS}

The paper proposes a systematic methodology to identify the core circuit organization of nonlinear behaviors frequently encountered in biological signalling. The circuit interpretation of the mathematical concept of organizing center emphasizes the role of feedback and monotonicity in organizing those behaviors. This insight dramatically simplifies the analysis of nonlocal attractors and might be useful for analysis and synthesis of  tunable yet robust nonlinear devices.




\section*{ACKNOWLEDGMENT}

Fernando Casta\~nos is gratefully acknowledged for insightful comments and suggestions during the visit of the first author to CINVESTAV, Mexico DF.

\addtolength{\textheight}{-0cm}   


\bibliographystyle{unsrt}
\bibliography{../../../../../../Dropbox/Orchestron_bibliography/orchestron}

\begin{thebibliography}{10}

\bibitem{Willems1976}
J.~C. Willems.
\newblock Realization of systems with internal passivity and symmetry
  constraints.
\newblock {\em Journal of the Franklin Institute}, 301(6):605--621, 1976.

\bibitem{Benvenuti2004}
L.~Benvenuti and L.~Farina.
\newblock A tutorial on the positive realization problem.
\newblock {\em IEEE Transactions on Automatic Control}, 49(5):651--664, 2004.

\bibitem{VanderSchaft1987}
A.~J. Van~der Schaft and P.~E. Crouch.
\newblock Hamiltonian and self-adjoint control systems.
\newblock {\em Systems \& control letters}, 8(4):289--295, 1987.

\bibitem{Golubitsky1985}
M.~Golubitsky and D.~G. Schaeffer.
\newblock {\em Singularities and Groups in Bifurcation Theory}, volume~51 of
  {\em Applied mathematical sciences}.
\newblock Springer-Verlag, New York, NY, 1985.

\bibitem{Franci2013b}
A.~Franci, G.~Drion, and R.~Sepulchre.
\newblock Modeling the modulation of neuronal bursting: a singularity theory
  approach.
\newblock {\em SIAM J Appl Dyn Syst}, 13(2):798--829, 2014.

\bibitem{Guckenheimer2002}
J.~Guckenheimer and P.~Holmes.
\newblock {\em Nonlinear oscillations, dynamical systems, and bifurcations of
  vector fields}, volume~42 of {\em Applied Mathematical Sciences}.
\newblock Springer, New-York, 7th edition, 2002.

\bibitem{Angeli2003}
D.~Angeli and E.~D. Sontag.
\newblock Monotone control systems.
\newblock {\em IEEE Transactions on Automatic Control}, 48(10):1684--1698,
  2003.

\bibitem{Hopfield1982}
J.~J. Hopfield.
\newblock Neural networks and physical systems with emergent collective
  computational abilities.
\newblock {\em Proceedings of the National Academy of Sciences},
  79(8):2554--2558, April 1982.

\bibitem{Hodgkin1952}
A.~Hodgkin and A.~Huxley.
\newblock A quantitative description of membrane current and its application to
  conduction and excitation in nerve.
\newblock {\em J. Physiol}, 117:500--544, 1952.

\bibitem{Griffith1968}
J.~S. Griffith.
\newblock Mathematics of cellular control processes {II}. positive feedback to
  one gene.
\newblock {\em J Theor Biol}, 20(2):209--16, 1968.

\bibitem{Griffith1968a}
J~S Griffith.
\newblock Mathematics of cellular control processes. {I}. negative feedback to
  one gene.
\newblock {\em J Theor Biol}, 20(2):202--8, 1968.

\bibitem{Angeli2004}
David Angeli, James~E. Ferrell, and Eduardo~D. Sontag.
\newblock Detection of multistability, bifurcations, and hysteresis in a large
  class of biological positive-feedback systems.
\newblock {\em Proc Natl Acad Sci U S A}, 101(7):1822--27, 2004.

\bibitem{Fenichel1971}
N.~Fenichel.
\newblock Persistence and smoothness of invariant manifolds for flows.
\newblock {\em Indiana University Mathematics Journal}, 21(3):193--226, 1971.

\bibitem{Krupa2001a}
M.~Krupa and P.~Szmolyan.
\newblock Extending slow manifolds near transcritical and pitchfork
  singularities.
\newblock {\em Nonlinearity}, 14:1473--1491, 2001.

\bibitem{Krupa2001b}
M.~Krupa and P.~Szmolyan.
\newblock Extending geometrical singular perturbation theory to nonhyperbolic
  points - folds and canards points in two dimensions.
\newblock {\em SIAM J. Math. Analysis}, 33(2):286--314, 2001.

\bibitem{Krupa2001c}
M.~Krupa and P.~Szmolyan.
\newblock Relaxation oscillation and canard explosion.
\newblock {\em J. Differential Equations}, 174(2):312--368, 2001.

\bibitem{Grasman1987}
J.~Grasman.
\newblock {\em Asymptotic methods for relaxation oscillations and
  applications}, volume~63.
\newblock Springer-Verlag New York, 1987.

\bibitem{Mishchenko1980}
E.~F. Mishchenko and N.~K. Rozov.
\newblock {\em Differential equations with small parameters and relaxation
  oscillations}, volume~13.
\newblock Plenum Publishing Corporation, 1980.

\bibitem{Angeli2008}
D.~Angeli and E.~D. Sontag.
\newblock Oscillations in i/o monotone systems under negative feedback.
\newblock {\em IEEE Transactions on Automatic Control}, 53(Special
  Issue):166--176, 2008.

\bibitem{Chow1994}
S.-N. Chow, C.~Li, and D.~Wang.
\newblock {\em Normal forms and bifurcation of planar vector fields}.
\newblock Cambridge University Press, 1994.

\bibitem{FitzHugh1961}
R.~FitzHugh.
\newblock Impulses and physiological states in theoretical models of nerve
  membrane.
\newblock {\em Biophysical J}, 1:445--466, 1961.

\bibitem{Rinzel1989}
J.~Rinzel and G.B. Ermentrout.
\newblock Analysis of neural excitability and oscillations.
\newblock In {\em Methods in neuronal modeling}, pages 135--169. MIT Press,
  1989.

\bibitem{Izhikevich2007}
E.~M. Izhikevich.
\newblock {\em Dynamical systems in neuroscience: the geometry of excitability
  and bursting}.
\newblock MIT Press, Cambridge, Mass., 2007.

\bibitem{Rinzel1987b}
J.~Rinzel.
\newblock A formal classification of bursting mechanisms in excitable systems.
\newblock In {\em Mathematical topics in population biology, morphogenesis and
  neurosciences}, pages 267--281. Springer, 1987.

\end{thebibliography}

\end{document}